\documentclass[a4paper]{amsart}

\usepackage[latin1]{inputenc}
\usepackage{amsthm, textcomp, amsmath, amsfonts}
\usepackage[pdftex]{graphicx}
\usepackage{url}
\usepackage{subfig}
\usepackage[a4paper,left=3cm,right=3cm]{geometry}
\allowdisplaybreaks[4]

\begin{document}

\title{Discrepancy of a convex set with zero curvature at one point} 
\author{B. Gariboldi}
\date{}

\begin{abstract}
Let $\Omega \subset \mathbb{R}^{d}$ be a convex body with everywhere positive curvature except at the origin and with
the boundary $\partial \Omega$ as the graph of the function $y=|x|^{\gamma}$ in a neighborhood of the origin with $\gamma \geq 2$. We consider the $L^{p}$ norm of the discrepancy with respect to translations and rotations of a dilated copy of the set $\Omega$.
\end{abstract}

\maketitle

\newtheorem{Theorem}{Theorem}
\newtheorem{Cor}[Theorem]{Corollary}
\newtheorem{lemma}[Theorem]{Lemma}
\newtheorem{definition}[Theorem]{Definition}

\theoremstyle{remark}
\newtheorem{remark}{Remark}

\section{Introduction}

Let $\Omega \subset \mathbb{R}^{d}$ be a convex body with everywhere positive curvature except at one single point, which we may assume to be the origin, with measure $|\Omega|$ and let
$$D_{R\Omega}=R^{d}|\Omega|-\mathrm{card}(\mathbb{Z}^{d} \cap R\Omega).$$
We assume that the boundary $\partial \Omega$ coincides with the graph of the function $y=|x|^{\gamma}$ in a neighborhood of the origin, with $\gamma\geq 2$. 

We consider the $L^{p}$ norm of the discrepancy with respect to translations and rotations of a dilated copy of the set $\Omega$:
$$\left\{ \int_{SO(d)}  \int_{T^{d}} |D_{R\Omega}(\sigma,t)|^{p}dx d\sigma \right\}^{1/p},$$
where
\[
D_{R\Omega}(\sigma,t)=R^{d}|\Omega|-\mathrm{card}\left(\mathbb{Z}^{d}\cap (R\sigma(\Omega)+t) \right).
\]

See \cite {BT} for a detailed description of several discrepancy problems in the planar case where the boundary coincides with the graph of the function $y=|x|^{\gamma}$ near the origin with $\gamma\geq 2$.

In order to put our results in an appropriate perspective, let us present a short non-exhaustive list of previous results. 
Kendall considered the mean square average of the discrepancy under translations and proved in \cite{Kendall} that if $\Omega$ is an oval in $\mathbb{R}^{d}$,
\[
\left\{\int_{T^{d}} |D_{R\Omega}(t)|^{2}dt \right\}^{1/2}\leq CR^{\frac{d-1}{2}}.
\]
Herz (see \cite{Herz}) and Hlawka (see \cite{Hlawka}) proved that
\[
\sup_{t \in T^{d}} \left\{|D_{R\Omega}(t)|\right\} \leq CR^{\frac{d(d-1)}{d+1}}.
\]
Huxley proved in \cite{Huxley} that if $\Omega$ is  a convex set in the plane with boundary with continuous
positive curvature, then
\[
\left\{ \int_{T^{2}} |D_{R\Omega}(t)|^{4}dt\right\}^{1/4}\leq CR^{1/2}\log^{1/4}(R).
\]
In \cite{cambridge} the authors extendend this result to higher dimension:
\[
\left\{ \int_{T^{d}} |D_{R\Omega}(t)|^{p}dt\right\}^{1/p} \leq 
\begin{cases}
CR^{\frac{d-1}{2}} & 1\leq p < 2d/(d-1),\\
CR^{\frac{d-1}{2}}\log^{\frac{d-1}{2d}}(R) & p=2d/(d-1),\\
CR^{\frac{d(d-1)}{d+1}\left(1-\frac{1}{p}\right)} & p> 2d/(d-1).
\end{cases}
\]
Mean square averages of the discrepancy under rotations of the domain have been considered by Iosevich in \cite{Iosevich}. Brandolini, Hofmann and Iosevich proved in \cite{BHI} that the average decay under rotations of the Fourier transform of the characteristic function of a convex set, without any further smoothness or curvature assumptions, satisfies the estimate
\[
\left\{\int_{SO(d)}|\widehat{\chi}_{\Omega}(\sigma\rho)|^{2}d\sigma \right\}^{1/2}\leq C(\mathrm{diameter}(\Omega))^{\frac{d-1}{2}}|\rho|^{-\frac{d+1}{2}},
\]
where $\rho\geq 1$.
Since the size of the discrepancy of a set is closely connected to the decay of its characteristic function, a corollary of this result in \cite{panorama} is that, without smoothness or curvature assumptions,
\[
\left\{\int_{SO(d)} \int_{T^{d}} |D_{R\Omega}(\sigma,t)|^{2}dtd\sigma \right\}^{1/2}\leq C(\mathrm{diameter}(\Omega))^{\frac{d-1}{2}}R^{\frac{d-1}{2}}.
\]
The results about the $L^{p}$ norm of the discrepancy with $p\neq 2$ are less complete. In particular, in \cite{polyhedra} the authors studied the $L^{p}$ norm of the discrepancy for rotated and translated polyhedra $\Omega$ in $\mathbb{R}^{d}$,
\[
\left\{\int_{SO(d)}\int_{T^{d}} |D_{R\Omega}(\sigma,t)|^{p} dtd\sigma \right\}^{1/p}\leq 
\begin{cases}
C\log^{d}(R) & p=1,\\
CR^{(d-1)\left(1-\frac{1}{p}\right)} & 1<p\leq +\infty.
\end{cases}
\]
In the same paper it was also observed that $L^{p}$ discrepancy of a polyhedron grows up with $p$. On the contrary, we see above that for certain domains with positive curvature there exists a range of indices $p$ where the $L^{p}$ discrepancy is of the same order of the $L^{2}$ discrepancy.


We now recall the definition of flat points:
\begin{definition}
Let $B$ a convex body in $\mathbb{R}^{d}$, let $z \in \partial B$ and let $\gamma>1$. We say that $z$ is an isolated flat point of order $\gamma$ if, in a neighborhood of $z$ and in a suitable Cartesian coordinate system with the origin in $z$, $\partial B$ is the graph of a function $\Phi$ such that
\begin{itemize}
\item[1.] $\Phi \in C^{\infty}(U\setminus\{ z\})$ where $U$ is a bounded open neighborhood of the origin $z$.
\item[2.] If, for every $x\in U\setminus\{ z\} $, $\mu_{1}(x), \cdots, \mu_{d-1}(x)$ are the eigenvalues of the Hessian matrix of $\Phi$,  then  for $j=1,\cdots,d-1$
\[
\inf_{x\in U\setminus\{ z\}} |x|^{2-\gamma}\mu_{j}(x)>0.
\]
\item[3.] For every multi-index $\alpha$,
\[
\sup_{x\in U\setminus\{ z\}} |x|^{|\alpha|-\gamma}\left\vert \frac{\partial^{|\alpha|}\Phi}{\partial x^{\alpha}}(x) \right\vert<+\infty.
\]
\end{itemize} 
\end{definition}
Notice that the convex set $\Omega$ we are considering is an example of convex body with a flat point at the origin.
In \cite{convex} the authors studied the discrepancy for convex bodies with a finite number of isolated flat points of order at most $\gamma$.  For $2<\gamma\leq d+1$ they proved that
\[
\left\{\int_{T^{d}} |D_{R\Omega}(t)|^{p}dt  \right\}^{1/p} \leq 
\begin{cases}
CR^{(d-1)\left(1-\frac{1}{\gamma} \right)} & 1\leq p\leq 2d/(d+1-\gamma),\\
CR^{\frac{d(d-1)}{d+1}\left(1-\frac{2}{\gamma p} \right)} & p> 2d/(d+1-\gamma);
\end{cases}
\]
for $\gamma > d+1$ they proved that
\[
\left\{\int_{T^{d}} |D_{R\Omega}(t)|^{p}dt  \right\}^{1/p} \leq CR^{(d-1)\left(1-\frac{1}{\gamma} \right)}.
\]
See also \cite{Colin, Guo1, Guo2, Guo3, ISS, Kendall, Kratzel, R1} and \cite{R2} for some results on the $L^{\infty}$ estimate of the discrepancy for specific classes of convex bodies with vanishing Gaussian curvature.

Different results are obtained with an average of the discrepancy with respect to translations and dilations:
\[
\left\{ \int_{\mathbb{R}}  \int_{T^{d}} |D_{r\Omega}(t)|^{p}dt dr \right\}^{1/p} \ \text{ or } \ \left\{\int_{T^{d}} \left(\int_{\mathbb{R}}   |D_{r\Omega}(t)|^{2}dr\right)^{p/2} dt \right\}^{1/p}.
\]
See \cite{normepure, normemiste} for some results for convex bodies with positive Gaussian curvature.

The present paper continues this line of research and it extends some results obtained in \cite{convex}. In particular, here we study
the average of the discrepancy function with respect to translations and rotations
of a dilated copy of a convex set $\Omega$ with everywhere positive curvature except at one single point and with
the boundary $\partial \Omega$ coinciding with the graph of the function $y=|x|^{\gamma}$ in a neighborhood of the origin, with $\gamma\geq 2$:
$$\left\{ \int_{SO(d)}  \int_{T^{d}} |D_{R\Omega}(\sigma,t)|^{p}dt d\sigma \right\}^{1/p}.$$

Notice that with the same techniques used in this paper, one can also estimate the $L^{s}(SO(d))\times L^{p}(T^{d})$ norm of the discrepancy with $s\neq p$:
\[
\left\{\int_{SO(d)} \left\{ \int_{T^{d}} |D_{R\Omega}(\sigma,t)|^{p}dt \right\}^{s/p} d\sigma \right\}^{1/s}.
\]

Our results are described by the following theorems.

\begin{Theorem} \label{teo1}
Let $\Omega$ be a convex body with everywhere positive curvature except at the origin and let the boundary $\partial \Omega$ be as the graph of the function $y=|x|^{\gamma}$ in a neighborhood of the origin with $\gamma \geq 2$. For $p>2$ one has
\begin{itemize}
\item[1)] If $2\leq\gamma\leq d+1$, 
\begin{align*}
&\left\{ \int_{SO(d)}  \int_{T^{d}} |D_{R\Omega}(\sigma,t)|^{p}dtd\sigma  \right\}^{1/p}\\
& \leq C
\begin{cases}
R^{(d-1)/2} & \textit{if } p<2d/(d-1),\\
R^{\frac{d(d-1)(p-2)}{d(p-2)+p}} & \textit{if } 2d/(d-1)\leq p< 2(\gamma-1)/(\gamma-2),\\
R^{\frac{d(d-1)}{d+\gamma-1}}(\log R)^{\frac{d-1}{p}} & \textit{if } p=2(\gamma-1)/(\gamma-2),\\
R^{\frac{d(d-1)}{d+1}\left(1-\frac{2(\gamma-1)}{p(d+\gamma-1)}\right)}(\log R)^{\frac{d-1}{p}} & \textit{if } p>2(\gamma-1)/(\gamma-2).
\end{cases} 
\end{align*}
In particular, if $\gamma=d+1$,
\begin{align*}
&\left\{ \int_{SO(d)}  \int_{T^{d}} |D_{R\Omega}(\sigma,t)|^{p}dtd\sigma  \right\}^{1/p}\\
& \leq C
\begin{cases}
R^{(d-1)/2} & \textit{if } p<2d/(d-1),\\
R^{\frac{d-1}{2}}(\log R)^{\frac{d-1}{p}} & \textit{if } p=2d/(d-1),\\
R^{\frac{d(d-1)}{d+1}\left(1-\frac{1}{p}\right)}(\log R)^{\frac{d-1}{p}} & \textit{if } p>2d/(d-1).
\end{cases} 
\end{align*}
\item[2)] If $\gamma>d+1$, 
$$\left\{ \int_{SO(d)}  \int_{T^{d}} |D_{R\Omega}(\sigma,t)|^{p}dtd\sigma  \right\}^{1/p} \leq C
\begin{cases}
R^{(d-1)/2} & \textit{if } p< 2(\gamma-1)/(\gamma-2),\\
R^{(d-1)/2}(\log R)^{\frac{d-1}{p}} & \textit{if } p=2(\gamma-1)/(\gamma-2),\\
R^{(d-1)\left(1-\frac{1}{p}\right) \left(1-\frac{1}{\gamma} \right)} & \textit{if } p> 2(\gamma-1)/(\gamma-2).
\end{cases} $$
\end{itemize}
\end{Theorem}

\begin{Cor}\label{teo2}
Let $\Omega$ be a convex body with everywhere positive curvature except at the origin and let the boundary $\partial \Omega$ be as the graph of the function $y=|x|^{\gamma}$ in a neighborhood of the origin with $\gamma \geq 2$. For $1\leq p\leq 2$ and for all $\gamma\geq 2$, one has
\[
\left\{ \int_{SO(d)}  \int_{T^{d}} |D_{R\Omega}(\sigma,t)|^{p}dtd\sigma  \right\}^{1/p} \leq CR^{(d-1)/2}
\]
\end{Cor}

Next theorem considers only rotations with $p=1$: the estimate  of the discrepancy averaged also over translations is better than this one, which coincides with the $L^{\infty}$ norm in the ball case. Notice that in this case, for $\gamma >d+1$, the estimate of the norm of the discrepancy with respect to rotations is better than the estimate with respect to translations (see \cite{convex}).

\begin{Theorem} \label{teo3}
Let $\Omega$ be a convex body with everywhere positive curvature except at the origin and let the boundary $\partial \Omega$ be as the graph of the function $y=|x|^{\gamma}$ in a neighborhood of the origin with $\gamma \geq 2$. For $p=1$ and for all $\gamma\geq 2$, one has
\[
\int_{SO(d)}  |D_{R\Omega}(\sigma)|d\sigma   \leq CR^{\frac{d(d-1)}{d+1}}.
\]
\end{Theorem}

It's easy to see that for every $p\geq 1$ the $L^{p}$ estimate of the discrepancy with respect only to rotations is always worse than the $L^{p}$ norm over rotations and translations. To prove it, it suffices to follow the proof of Theorem \ref{teo3} replacing the $L^{1}$ norm with the $L^{p}$ one.

We can also estimate the discrepancy from below, in every dimension. 

\begin{Theorem} \label{teo4}
Let $\Omega$ be a convex body with everywhere positive curvature except at the origin and let the boundary $\partial \Omega$ be as the graph of the function $y=|x|^{\gamma}$ in a neighborhood of the origin with $\gamma \geq 2$. Then, for all $\gamma \geq 2$
\[
\left\{\int_{SO(d)}\int_{T^{d}}|D_{R\Omega}(\sigma, t)|^{p}dt d\sigma\right\}^{1/p}\geq cR^{\frac{d-1}{2}}.
\]
\end{Theorem}

\begin{figure}[!h] 
\centering
 \includegraphics[width=16 cm, height=10 cm]{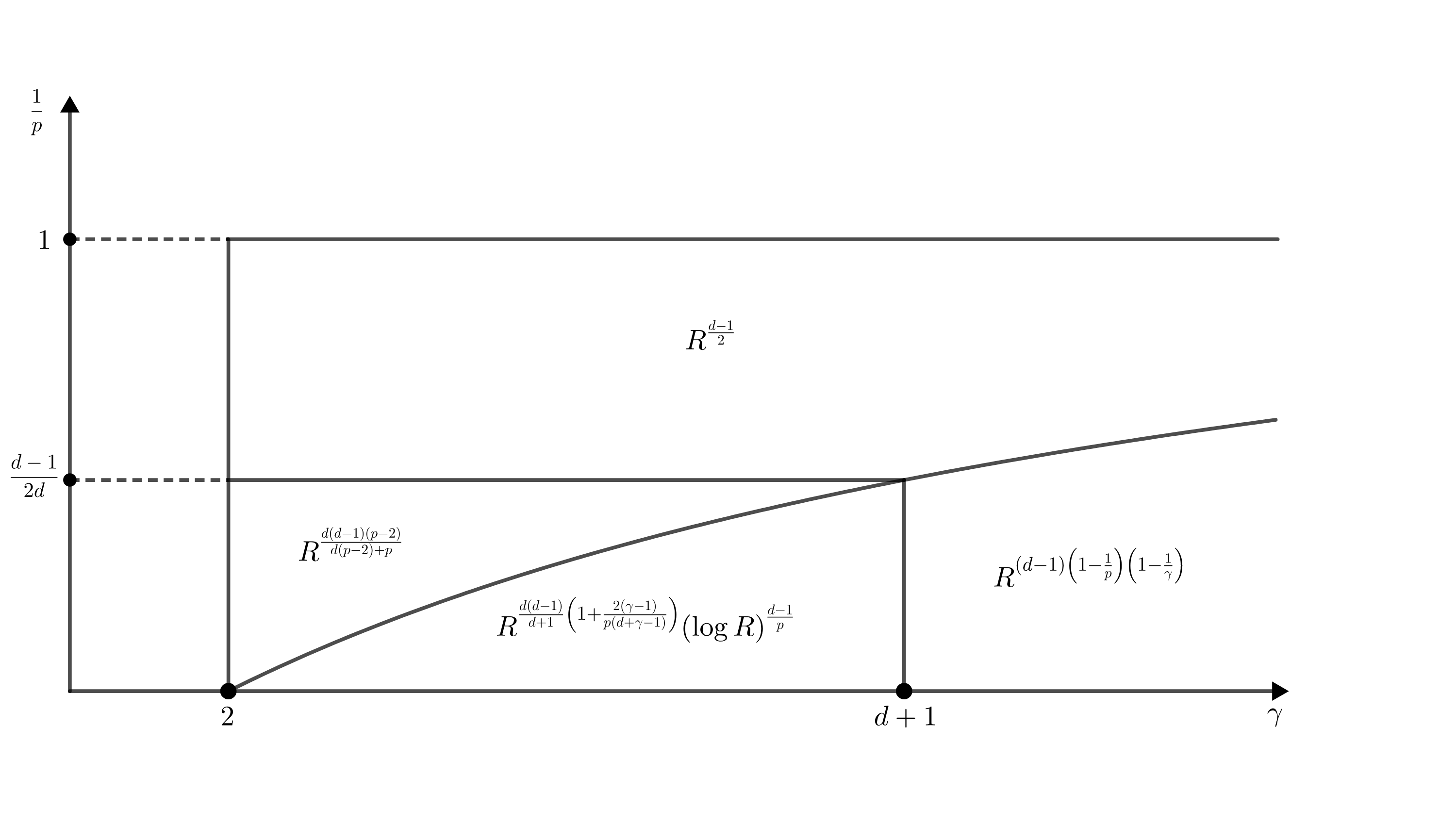} 
\caption{}\label{fig1}
 \end{figure}

Notice that, when $\gamma< d+1$ for $p<2d/(d-1)$ and when $\gamma=d+1$ for every $p$,  estimates remember the ones in \cite{cambridge} for the case of positive Gaussian curvature and are independent of $\gamma$: rotations do not influence so much the results. Notice also that if $p=+\infty$ the estimates for $L^{\infty}$ discrepancy of $\Omega$ when $\gamma \leq d+1$ match Landau's estimates for the $L^{\infty}$ discrepancy of the ball (see \cite{Landau}) and also when $\gamma > d+1$ correspond to the results in \cite{convex}, where the authors do not consider rotations.

Figure \ref{fig1} summarizes the estimates for the discrepancy.
\newline

I would like to thank Giancarlo Travaglini
for his valuable suggestions during the preparation of this paper.

\section{Preliminary lemmas}

D. G. Kendall observed that the discrepancy is a periodic function of the
translation, and it has a Fourier expansion with coefficients that are a
sampling of the Fourier transform of $\Omega$,
\[
\widehat{\chi}_{\Omega}\left(  \xi\right)  =%
{\displaystyle\int_{\Omega}}
\exp\left(  -2\pi i\xi x\right)  dx.
\]

\begin{lemma}
\label{Fourier} The number of integer points in $R\Omega-t$, a translated by a
vector $t\in\mathbb{R}^{d}$\ and dilated by a factor $R\geq 1$\ of a domain
$\Omega$\ in the $d$ dimensional Euclidean space, is a periodic function of the
translation with Fourier expansion
\[
\sum_{k\in\mathbb{Z}^{d}}
\chi_{R\Omega-t}(k)=\sum_{n\in\mathbb{Z}^{d}}
R^{d}\widehat{\chi}_{\Omega}\left(  Rn\right)  \exp(2\pi int).
\]
In particular, the discrepancy is
\[
{\displaystyle\sum_{n\in\mathbb{Z}^{d}\setminus\left\{  0\right\}  }}
R^{d}\widehat{\chi}_{\Omega}\left(  Rn\right)  \exp(2\pi int).
\]

\end{lemma}

\proof
Indeed
\begin{align*}
\sum_{k\in\mathbb{Z}^{d}} \chi_{R\Omega-t}(k) & = \sum_{n \in \mathbb{Z}^{d}} \left(\int_{T^{d}} \sum_{k\in \mathbb{Z}^{d}} \chi_{R\Omega-y}(k)\exp(-2\pi iny)dy \right)\exp(2\pi int)\\
& \sum_{n \in \mathbb{Z}^{d}} \left(\int_{T^{d}} \chi_{R\Omega}(y)\exp(-2\pi iny)dy \right)\exp(2\pi int)\\
& \sum_{n \in \mathbb{Z}^{d}} R^{d}\widehat{\chi}_{\Omega}(Rn)\exp(2\pi int).
\end{align*}
\endproof

Therefore
\begin{align*}
D_{R\Omega}(\sigma,t)&=R^{d}|\Omega|-\mathrm{card}\left(\mathbb{Z}^{d}\cap (R\sigma(\Omega)+t) \right)\\
&= R^{d}|\Omega| -\sum_{n\in \mathbb{Z}^{d}} \widehat{\chi}_{R\Omega}(\sigma(n))\exp(2\pi int)\\
&= R^{d}|\Omega| -R^{d}\sum_{n\in \mathbb{Z}^{d}} \widehat{\chi}_{\Omega}(R\sigma(n))\exp(2\pi int). 
\end{align*}

\begin{remark}\label{r1}
We emphasize that the Fourier expansion of the discrepancy converges at least
in $L^{2}\left(  \mathbb{T}^{d}\right)  $, but we are not claiming that it
converges pointwise. Indeed, the discrepancy is discontinuous, hence the
associated Fourier expansion does not converge absolutely or uniformly. To
overcome this problem, one can introduce a mollified discrepancy (see \cite{Giancarlo}). If the
domain $\Omega$\ is convex and it contains the origin, then there exists
$\varepsilon>0$\ such that if $\varphi\left(  x\right)  $\ is a non negative
smooth radial function with support in $\left\{  \left\vert x\right\vert
\leq\varepsilon\right\}  $\ and with integral 1, and if
$R\geq1$, then by \cite[Lemma 29]{panorama}
\[
\varepsilon^{-d}\varphi(\varepsilon^{-1}\cdot)\ast\chi_{(R-\varepsilon)\Omega}(x)
\leq\chi_{R\Omega}(x)\leq
\varepsilon^{-d}\varphi(\varepsilon^{-1}\cdot)\ast\chi_{(R+\varepsilon)\Omega}(x)
\]
and therefore

\begin{gather*}
\left\vert \Omega\right\vert \left(  \left(  R-\varepsilon\right)  ^{d}%
-R^{d}\right)  +\left(  R-\varepsilon\right)  ^{d}%
{\displaystyle\sum_{n\in\mathbb{Z}^{d}\setminus\left\{  0\right\}  }}
\widehat{\varphi}\left(  \varepsilon \sigma(n)\right)  \widehat{\chi}_{\Omega}\left(
\left(  R-\varepsilon\right)  \sigma(n)\right)  \exp\left(  2\pi int\right)\\
\leq D_{R\Omega}(\sigma,t)\\
\leq\left\vert \Omega\right\vert \left(  \left(  R+\varepsilon\right)  ^{d}%
-R^{d}\right)  +\left(  R+\varepsilon\right)  ^{d}%
{\displaystyle\sum_{n\in\mathbb{Z}^{d}\setminus\left\{  0\right\}  }}
\widehat{\varphi}\left(  \varepsilon \sigma(n)\right)  \widehat{\chi}_{\Omega}\left(
\left(  R+\varepsilon\right)  \sigma(n)\right)  \exp\left(  2\pi int\right) .
\end{gather*}
One has $\left\vert \left(  R+\varepsilon\right)  ^{d}-R^{d}\right\vert \leq
CR^{d-1}\varepsilon$, and one can work with the mollified discrepancy defined as
\[
D_{R\Omega}^{\varepsilon}(\sigma,t)=R  ^{d}%
{\displaystyle\sum_{n\in\mathbb{Z}^{d}\setminus\left\{  0\right\}  }}
\widehat{\varphi}\left(  \varepsilon \sigma(n)\right)  \widehat{\chi}_{\Omega}\left(
R \sigma(n)\right)  \exp\left(  2\pi int\right)  .
\]
In particular 
\[
\vert D_{R\Omega}(\sigma,t)\vert \leq CR^{d-1}\varepsilon+R^{d}\sum_{n\in\mathbb{Z}^{d}\setminus\left\{  0\right\}  }
\widehat{\varphi}\left(  \varepsilon \sigma(n)\right)  \widehat{\chi}_{\Omega}\left(
R \sigma(n)\right)  \exp\left(  2\pi int\right) .
\]
Observe that since $\left\vert \widehat{\varphi
}\left(  \zeta\right)  \right\vert \leq C\left(  1+\left\vert \zeta\right\vert
\right)  ^{-\gamma}$ for every $\gamma>0$, the mollified Fourier
expansion has no problems of convergence.
\end{remark}

Another useful result comes from the decay of the Fourier transform in \cite{convex, Greenleaf, Rigoli}, with $(\xi',\xi_{d})\in \mathbb{R}^{d-1}\times\mathbb{R}$
\begin{equation} \label{fourier}
|\widehat{\chi}_{\Omega}(\xi',\xi_{d})|\leq \begin{cases}
c|\xi_{d}|^{-1-\frac{d-1}{\gamma}},\\
c|\xi'|^{-(d-1)\frac{\gamma-2}{2(\gamma-1)}}|\xi_{d}|^{-\frac{d-1}{2(\gamma-1)}-1},\\
c|\xi'|^{-\frac{d+1}{2}}.
\end{cases}
\end{equation}

\begin{Theorem}
Let $\Omega$ be a convex body with everywhere positive curvature except at one single point and with
the boundary $\partial \Omega$ coinciding with the graph of the function $y=|x|^{\gamma}$ in a neighborhood of the origin, with $\gamma\geq 2$. Then, if $\rho>1$ 
\begin{equation}\label{rot}
\left\{ \int_{SO(d)} |\widehat{\chi}_{\Omega}(\rho \omega)|^{p} d\omega \right\}^{1/p} \leq 
\begin{cases}
C\rho^{-(d+1)/2} & \text{if } p< 2(\gamma-1)/(\gamma-2),\\
C\rho^{-(d+1)/2}\log^{(\gamma-2)(d-1)/2(\gamma-1)}(\rho) & \text{if } p= 2(\gamma-1)/(\gamma-2),\\
C\rho^{-(d-1)\left(\frac{1}{p}+\frac{1}{\gamma}-\frac{1}{p\gamma} \right)-1}  & \text{if } p> 2(\gamma-1)/(\gamma-2).\\
\end{cases}
\end{equation}
\end{Theorem}

\proof Let $(\xi',\xi_{d})=(\rho \omega'\sin\theta, \rho \cos\theta)$ with $\omega' \in SO(d-1)$ and $0\leq \theta \leq \pi$. When $\varepsilon<\theta<\pi-\varepsilon$ one has the classical estimate 
\[
|\widehat{\chi}_{\Omega}(\rho \omega'\sin\theta, \rho \cos\theta)|\leq \rho^{-\frac{d+1}{2}}.
\]
Therefore for $p>2(\gamma-1)/(\gamma-2)$ one has
\begin{align*}
&\int_{0}^{\pi}\int_{SO(d-1)} |\widehat{\chi}_{\Omega}(\rho \omega'\sin\theta, \rho \cos\theta)|^{p}\sin^{d-2}\theta d\omega'd\theta\\
&\leq c\rho^{-\frac{d+1}{2}p}+c\int_{0}^{\varepsilon}\left\{\min\left(\rho^{-\frac{d-1}{\gamma}-1}(\cos\theta)^{-\frac{d-1}{\gamma}-1},\rho^{-\frac{d+1}{2}}(\sin\theta)^{-(d-1)\frac{\gamma-2}{2(\gamma-1)}}(\cos\theta)^{-\frac{d-1}{2(\gamma-1)}-1} \right) \right\}^{p}\sin^{d-2}\theta d\theta\\
&\leq c\rho^{-\frac{d+1}{2}p}+c\int_{0}^{\varepsilon} \left\{\min\left(\rho^{-\frac{d-1}{\gamma}-1},\rho^{-\frac{d+1}{2}}\theta^{-(d-1)\frac{\gamma-2}{2(\gamma-1)}}\right) \right\}^{p}\theta^{d-2}d\theta\\
&= c\rho^{-\frac{d+1}{2}p}+c\int_{0}^{\rho^{-1+1/\gamma}} \rho^{\left(-\frac{d-1}{\gamma}-1\right)p}\theta^{d-2}d\theta+c\int_{\rho^{-1+1/\gamma}}^{\varepsilon} \rho^{-\frac{d+1}{2}p}\theta^{-(d-1)\frac{\gamma-2}{2(\gamma-1)}p}\theta^{d-2}d\theta\\
&= c\rho^{-\frac{d+1}{2}p}+c\rho^{-(d-1)\frac{\gamma+p-1}{\gamma}p-p}+\rho^{-(d-1)\frac{\gamma+p-1}{\gamma}p-p}\leq \rho^{-(d-1)\frac{\gamma+p-1}{\gamma}p-p}.
\end{align*}
In the same way, one can obtain the other two results. \endproof

\section{Estimates of the $L^{\infty}$ norm}

If $\Omega$ is a convex set with positive Gaussian curvature, from \cite{cambridge} and \cite{Hlawka} we know that
\[
\sup_{t\in T^{d}} |D_{R\Omega}(t)| \leq CR^{d-2+\frac{2}{d+1}}.
\]
In the case we consider a convex set $\Omega$ with a finite number of isolated flat points of order at most $\gamma$ (see \cite{convex, Colin}), this estimate holds for $2<\gamma\leq d+1$. Instead, for $\gamma >d+1$, we have
\[
\sup_{t\in T^{d}} |D_{R\Omega}(t)| \leq CR^{(d-1)\left(1-\frac{1}{\gamma} \right)}.
\]
In this section, we consider the $L^{\infty}$ norm of the discrepancy with respect to translations and rotations of a convex set $\Omega$ with everywhere positive curvature except at the origin and with
the boundary $\partial \Omega$ as the graph of the function $y=|x|^{\gamma}$ in a neighborhood of the origin.
Notice from the following result that rotations do not influence the estimate of the $L^{\infty}$ norm.

\begin{Theorem}
Let $\Omega \subset \mathbb{R}^{d}$ be a convex body with everywhere positive curvature except at the origin and with
the boundary $\partial \Omega$ as the graph of the function $y=|x|^{\gamma}$ in a neighborhood of the origin.
Let $\gamma\geq 2$. Then
\[
\sup_{t \in T^{d}}|D_{R\Omega}(\sigma,t)|\leq \begin{cases}
C_{\gamma,d} R^{d-2+\frac{2}{d+1}} & \gamma\leq d+1\\
C_{\gamma,d} R^{(d-1)\left(1-\frac{1}{\gamma} \right)} & \gamma> d+1.
\end{cases}
\]
\end{Theorem}

\proof One has
\begin{equation}\label{dis}
|D_{R\Omega}(\sigma,t)|\leq R^{d-1}\varepsilon + R^{d}\sum_{n\neq 0} |\widehat{\chi}_{\Omega}(R\sigma n)|\frac{1}{1+|\varepsilon n|^{K}}.
\end{equation}
Let $\xi=(\xi',\xi_{d})$ with $\xi' \in \mathbb{R}^{d-1}$ and $\xi_{d}\in \mathbb{R}$. One has
$$(\xi',\xi_{d})=(\rho \omega' \sin \theta, \rho \cos \theta)$$
with $\rho >0$, $\omega' \in S^{d-2}$ and $0\leq \theta\leq \pi$. 
We want to replace the series in (\ref{dis}) with the following integral: 
$$\int_{|\xi|\geq 1/2} |\widehat{\chi}_{\Omega}(Rn)|\frac{1}{1+|\varepsilon n|^{K}}d\xi.$$
This means that one has to consider the set of points near the origin and the set of the line $\{n_{d}=0\}$ in a different way. Hence, remembering that $\partial \Omega$ has a strictly positive curvature out of a neighborhood of the origin, for \ref{Fourier} one has
\begin{align*}
 |D_{R\Omega}(\sigma,t)|& \leq R^{d-1}\varepsilon +R^{d}\max_{1\leq |n|\leq c}|\widehat{\chi}_{\Omega}(R n)|+R^{d}\sum_{n_{d}=1}^{+\infty}|\widehat{\chi}_{\Omega}(R(0,\cdots,n_{d}))|\\
& + R^{d}\int_{1/2}^{+\infty}\int_{0}^{\pi}\int_{S^{d-2}} |\widehat{\chi}_{\Omega}(R\rho(\omega'\sin\theta,\cos\theta))|d\omega'\sin^{d-2}(\theta)d\theta \frac{1}{1+(\varepsilon \rho)^{K}}\rho^{d-1}d\rho \\
& \leq R^{d-1}\varepsilon+R^{(d-1)\left(1-\frac{1}{\gamma} \right)}+R^{(d-1)\left(1-\frac{1}{\gamma} \right)}\sum_{n_{d}=1}^{+\infty}n_{d}^{-\frac{d-1}{\gamma}-1}+R^{\frac{d(d-1)}{d+1}}\\
&+ R^{d}\int_{1/2}^{+\infty}\int_{0}^{\pi/4} \min\left((R\rho)^{-\frac{d-1}{\gamma}-1}, (R\rho\sin\theta)^{-\frac{(d-1)(\gamma-2)}{2(\gamma-1)}}(R\rho\cos\theta)^{-\frac{d-1}{2(\gamma-1)}-1}\right)\\
&\times \sin^{d-2}\theta d\theta \frac{1}{1+(\varepsilon \rho)^{K}}\rho^{d-1}d\rho.
\end{align*}
If $0\leq \theta\leq \pi/4$, one has
$$(R\rho)^{-\frac{d-1}{\gamma}-1}= (R\rho\sin\theta)^{-\frac{(d-1)(\gamma-2)}{2(\gamma-1)}}(R\rho\cos\theta)^{-\frac{d-1}{2(\gamma-1)}-1}$$
$$\Rightarrow \theta \approx (R\rho)^{\frac{1}{\gamma}-1}.$$
Therefore one has to consider the following integral:
\begin{align*}
 & R^{d}\int_{1/2}^{+\infty}\int_{0}^{c(R\rho)^{\frac{1}{\gamma}-1}} \theta^{d-2}d\theta (R\rho)^{-\frac{d-1}{\gamma}-1}\frac{1}{1+(\varepsilon \rho)^{K}}\rho^{d-1}d\rho\\ 
&  + R^{d}\int_{1/2}^{+\infty}\int_{c(R\rho)^{\frac{1}{\gamma}-1}}^{\pi/4} (R\rho\theta)^{-\frac{(d-1)(\gamma-2)}{2(\gamma-1)}}(R\rho)^{-\frac{d-1}{2(\gamma-1)}-1}\theta^{d-2}d\theta\\
& \times \frac{1}{1+(\varepsilon \rho)^{K}}\rho^{d-1}d\rho\\
& := A+B.
\end{align*}
One has
\begin{align*}
A &  =R^{d}\int_{1/2}^{+\infty}\int_{0}^{c\left(  R\rho\right)  ^{\frac{1}{\gamma}-1}}\theta^{d-2}\ d\theta\ \left(  R\rho\right)  ^{-\frac{d-1}{\gamma}-1}\frac{1}{\left(  1+\varepsilon\rho\right)  ^{K}}\rho
^{d-1}\ d\rho\\
&  =R^{\left(  d-1\right)  \left(  1-\frac{1}{\gamma}\right)  }\int
_{1/2}^{+\infty}\int_{0}^{c\left(  R\rho\right)  ^{\frac{1}{\gamma}-1}}%
\theta^{d-2}\ d\theta\ \rho^{-\frac{d-1}{\gamma}+d-2}\frac{1}{\left(
1+\varepsilon\rho\right)  ^{K}}\ d\rho\\
&  \leq R^{\left(  d-1\right)  \left(  1-\frac{1}{\gamma}\right)  }%
\int_{1/2}^{+\infty}\ \left(  R\rho\right)  ^{-\frac{\left(  \gamma-1\right)
\left(  d-1\right)  }{\gamma}}\ \rho^{-\frac{d-1}{\gamma}+d-2}\frac{1}{\left(
1+\varepsilon\rho\right)  ^{K}}\ d\rho\\
&  =\int_{1/2}^{+\infty}\rho^{-1}\frac{1}{\left(  1+\varepsilon\rho\right)
^{K}}\ d\rho\leq\int_{\varepsilon/2}^{+\infty}t^{-1}\frac{1}{\left(
1+t\right)  ^{K}}\ dt\leq\log\left(  \frac{1}{\varepsilon}\right)  \ .
\end{align*}
and
\begin{align*}
B &  =R^{d}\int_{1/2}^{+\infty}\int_{c\left(  R\rho\right)  ^{\frac{1}{\gamma
}-1}}^{\pi/4}\left(  R\rho\theta\right)  ^{-\frac{\left(  d-1\right)  \left(
\gamma-2\right)  }{2\left(  \gamma-1\right)  }}\left(  R\rho\right)
^{-\frac{d-1}{2\left(  \gamma-1\right)  }-1}\theta^{d-2}\ d\theta \frac{1}{\left(  1+\varepsilon\rho\right)  ^{K}}\rho^{d-1}\ d\rho\\
&  =R^{\frac{d-1}{2}}\int_{1/2}^{+\infty}\int_{c\left(  R\rho\right)
^{\frac{1}{\gamma}-1}}^{\pi/4}\theta^{\frac{d\gamma+2-3\gamma}{2\left(
\gamma-1\right)  }}\ d\theta\ \rho^{\frac{d-3}{2}}\frac{1}{\left(
1+\varepsilon\rho\right)  ^{K}}\ d\rho\\
&  \leq R^{\frac{d-1}{2}}\int_{1/2}^{+\infty}\ \rho^{\frac{d-3}{2}}%
\frac{1}{\left(  1+\varepsilon\rho\right)  ^{K}}\ d\rho=R^{\frac{d-1}{2}}%
\int_{\varepsilon}^{+\infty}t^{\frac{d-3}{2}}\varepsilon^{-\frac{d-3}{2}}%
\frac{1}{\varepsilon\left(  1+t\right)  ^{K}}\ dx\\
&  \leq R^{\frac{d-1}{2}}\varepsilon^{-\frac{d-1}{2}}\int_{0}^{+\infty
}t^{\frac{d-3}{2}}\frac{1}{\varepsilon\left(  1+t\right)  ^{K}}\ dt\leq
R^{\frac{d-1}{2}}\varepsilon^{-\frac{d-1}{2}}\ ,
\end{align*}
because $\theta^{\frac{d\gamma+2-3\gamma}{2\left(  \gamma-1\right)  }}$ is convergent in $0^{+}$. Therefore one has
$$\left\vert D_{R\Omega}(\sigma,t)\right\vert \leq R^{d-1}\varepsilon+R^{\left(  d-1\right)
\left(  1-\frac{1}{\gamma}\right)  }+R^{\frac{d\left(  d-1\right)  }{d+1}%
}+\log\left(  \frac{1}{\varepsilon}\right)  +R^{\frac{d-1}{2}}\varepsilon
^{-\frac{d-1}{2}}\ ,$$
and, if $R^{d-1}\varepsilon=R^{\frac{d-1}{2}}\varepsilon^{-\frac{d-1}{2}}$, we have $\varepsilon=R^{-\frac{d-1}{d+1}}$ and
\begin{align*}
\left\vert D_{R\Omega}(\sigma,t)\right\vert &\leq R^{\left(  d-1\right)  \left(  1-\frac
{1}{\gamma}\right)  }+\log\left(  R\right)  +R^{\frac{d\left(  d-1\right)
}{d+1}}\\
&\leq \begin{cases}
 R^{\frac{d\left(  d-1\right)  }{d+1}} = C_{\gamma,d} R^{d-2+\frac{2}{d+1}} & \gamma\leq d+1\\
R^{\left(  d-1\right)  \left(  1-\frac{1}{\gamma}\right)  } = C_{\gamma,d}R^{\left(  d-1\right)  \left(  1-\frac
{1}{\gamma}\right)  } & \gamma> d+1.
\end{cases} 
\end{align*}
\endproof

\section{Proofs of the results}

Now let $q$ be such that $1/p+1/q=1$ with $p>2$. By the definition of discrepancy and the Hausdorff-Young inequality, one has
\begin{align*}
\left\{ \int_{SO(d)}  \int_{T^{d}} |D_{R\Omega}(\sigma,t)|^{p}dtd\sigma  \right\}^{1/p} & =
\left\{ \int_{SO(d)} \left\{ \left(\int_{T^{d}} |D_{R\Omega}(\sigma,t)|^{p}dt\right)^{1/p} \right\}^{p} d\sigma \right\}^{1/p}\\
&= \left\{ \int_{SO(d)} \left\{ \left(R^{d} \int_{T^{d}} \left| \sum_{n\neq 0} \widehat{\chi}_{\Omega}(R\sigma(n))e^{2\pi int}\right|^{p}dt\right)^{1/p}\right\}^{p}d\sigma\right\}^{1/p}\\
&\leq R^{d}\left\{ \int_{SO(d)} \left( \sum_{n\neq 0}\left| \widehat{\chi}_{\Omega}(R\sigma(n))\right|^{q}\right)^{p/q} d\sigma\right\}^{1/p}.
\end{align*}
Using the $L^{p/q}(SO(d))$-norm (note that $p/q>1$) and the Minkowski inequality, one has
\begin{align*}
R^{d}\left\{ \int_{SO(d)} \left( \sum_{n\neq 0}\left| \widehat{\chi}_{\Omega}(R\sigma(n))\right|^{q}\right)^{p/q} d\sigma\right\}^{q/p\cdot 1/q}
&=R^{d} \left\{ \left\| \sum_{n\neq 0} \left| \widehat{\chi}_{\Omega}(R\sigma(n))\right|^{q}\right\|_{L^{p/q}(SO(d))}\right\}^{1/q}\\
&\leq R^{d} \left\{ \sum_{n\neq 0}\left\|  \left| \widehat{\chi}_{\Omega}(R\sigma(n))\right|^{q}\right\|_{L^{p/q}(SO(d))}\right\}^{1/q}\\
&=R^{d}\left\{ \sum_{n\neq 0} \left\{ \int_{SO(d)} \left| \widehat{\chi}_{\Omega}(R\sigma(n))\right|^{p}d\sigma\right\}^{q/p}\right\}^{1/q}\\
&=R^{d} \left\{ \sum_{n\neq 0} \left\| \widehat{\chi}_{\Omega}(R\sigma(n))\right\|_{L^{p}(SO(d))}^{q}\right\}^{1/q}.
\end{align*}
In this way, to estimate the discrepancy, one can use the estimate for the Fourier transform of the characteristic function with respect to rotations (see (\ref{rot})).

\proof[Proof of Theorem \ref{teo1}] We start proving 1). First let $p<2(\gamma-1)/(\gamma-2)$. Because of (\ref{rot}), for $p<2d/(d-1)$, one has
$$R^{d} \left\{ \sum_{n\neq 0} \left\| \widehat{\chi}_{\Omega}\right\|_{L^{p}(SO(d))}^{q}\right\}^{1/q}
 \leq R^{(d-1)/2}\left\{\sum_{n\neq 0} |n|^{-\frac{d+1}{2}q} \right\}^{1/q}\\
\leq R^{(d-1)/2}.$$
Notice that  $2d/(d-1) < 2(\gamma-1)/(\gamma-2)$ if $\gamma<d+1$. Therefore, one has to study the case $2d/(d-1)\leq p<2(\gamma-1)/(\gamma-2)$. One can use  the mollified discrepancy with a cut-off function $\varphi$:
$$|D_{R\Omega}(\sigma,t)| \leq CR^{d-1}\varepsilon +R^{d}\left|\sum_{n\neq 0}\widehat{\varphi}(\varepsilon \sigma(n))\widehat{\chi}_{\Omega}(R\sigma(n)) \right|.$$
For $p>2d/(d-1)$,  with $\widehat{\varphi}(\xi)\leq (1+|\xi|)^{-K}$, one has
\begin{align*}
R^{d}\left\{ \int_{SO(d)}  \left|\sum_{n\neq 0}\widehat{\varphi}(\varepsilon \sigma(n))\widehat{\chi}_{\Omega}(R\sigma(n)) \right|^{p}d\sigma\right\}^{1/p}
&\leq R^{d}\left\{ \sum_{n\neq 0} (1+|\varepsilon n|)^{-K} \left\|\widehat{\chi}_{\Omega}(R\sigma(n)) \right\|^{q}_{L^{p}(SO(d))} \right\}^{1/q}\\
&\leq R^{d} \left\{\int_{|\xi|\geq 1} (1+\varepsilon |\xi|)^{-K}\left(R|\xi| \right)^{-\frac{d+1}{2}q} d\xi\right\}^{1/q}.
\end{align*}
Using polar coordinates and then the change of variables $\varepsilon \rho=y$, one obtains
\begin{align*}
&R^{(d-1)/2} \left\{\int_{1}^{+\infty} (1+\varepsilon \rho)^{-K}\rho^{-\frac{d+1}{2}q+d-1} d\rho\right\}^{1/q}\\
&\leq R^{(d-1)/2}\varepsilon^{(d+1)/2 - d/q} \left\{\int_{\varepsilon}^{+\infty} (1+y)^{-K}y^{-\frac{d+1}{2}q+d-1} dy\right\}^{1/q}\\
&\leq R^{(d-1)/2}\varepsilon^{(d+1)/2 - d/q}.
\end{align*}
Hence
$$\left\{ \int_{SO(d)}\int_{T^{d}} |D_{R\Omega}(\sigma,t)|^{p}dt d\sigma\right\}^{1/p}\leq C\left(R^{d-1}\varepsilon+R^{(d-1)/2}\varepsilon^{(d+1)/2 - d/q}\right)$$
and choosing $\varepsilon=R^{\frac{(d-1)p}{2d-p(d+1)}}$, for $\gamma<d+1$ and $2d/(d-1)<p<2(\gamma-1)/(\gamma-2)$ one has
$$\left\{ \int_{SO(d)}\int_{T^{d}} |D_{R\Omega}(\sigma,t)|^{p}dt d\sigma\right\}^{1/p}\leq CR^{\frac{d(d-1)(p-2)}{d(p-2)+p}}.$$

If $p=2(\gamma-1)/(\gamma-2)$, using the cut-off function and the corresponding estimate in (\ref{rot}), one has
\begin{align*}
&R^{d}\left\{ \sum_{n\neq 0} (1+|\varepsilon n|)^{-K} \left\|\widehat{\chi}_{\Omega}(R\sigma(n)) \right\|^{q}_{L^{p}(SO(d))} \right\}^{1/q}\\
&\leq R^{d} \left\{\int_{|\xi|\geq 1} (1+\varepsilon |\xi|)^{-K}\left(R|\xi| \right)^{-\frac{d+1}{2}q}\left(\log(R|\xi|)\right)^{\frac{(\gamma-2)(d-1)}{2(\gamma-1)}q} d\xi\right\}^{1/q}\\
&\leq R^{(d-1)/2} \left\{\int_{1}^{+\infty} (1+\varepsilon \rho)^{-K}\rho^{-\frac{d+1}{2}q+d-1}\left(\log(R\rho)\right)^{\frac{(\gamma-2)(d-1)}{2(\gamma-1)}q} d\rho\right\}^{1/q}\\
&\leq R^{(d-1)/2}\varepsilon^{(d+1)/2 - d/q} \left\{\int_{\varepsilon}^{+\infty} (1+y)^{-K}y^{-\frac{d+1}{2}q+d-1}\left(\log(Ry/\varepsilon)\right)^{\frac{(\gamma-2)(d-1)}{2(\gamma-1)}q} dy\right\}^{1/q}.
\end{align*}
Hence
\begin{align*}
&\int_{\varepsilon}^{+\infty} (1+y)^{-K}y^{-\frac{d+1}{2}q+d-1}\left(\log(Ry/\varepsilon)\right)^{\frac{(\gamma-2)(d-1)}{2(\gamma-1)}q} dy\\
& \leq \int_{\varepsilon}^{+\infty} (1+y)^{-K}y^{-\frac{d+1}{2}q+d-1}\left(\log(R/\varepsilon)\right)^{\frac{(\gamma-2)(d-1)}{2(\gamma-1)}q} dy\\
& + \int_{\varepsilon}^{+\infty} (1+y)^{-K}y^{-\frac{d+1}{2}q+d-1}\left(\log(y)\right)^{\frac{(\gamma-2)(d-1)}{2(\gamma-1)}q} dy\\
& := A+B.
\end{align*}
It's easy to see that 
\[
A\leq \log(R/\varepsilon)^{\frac{(\gamma-2)(d-1)}{2(\gamma-1)}q}
\]
and that $B$ converges for $\gamma<d+1$.
Therefore
\[
R^{d}\left\{ \sum_{n\neq 0} (1+|\varepsilon n|)^{-K} \left\|\widehat{\chi}_{\Omega}(R\sigma(n)) \right\|^{q}_{L^{p}(SO(d))} \right\}^{1/q}\leq CR^{(d-1)/2}\varepsilon^{(d+1)/2 - d/q}\log(R/\varepsilon)^{\frac{(d-1)}{p}}
\]
and choosing $\varepsilon=R^{\frac{(d-1)p}{2d-p(d+1)}}$, for $\gamma<d+1$ and $p=2(\gamma-1)/(\gamma-2)$ one obtains
$$\left\{ \int_{SO(d)}\int_{T^{d}} |D_{R\Omega}(\sigma,t)|^{p}dt d\sigma\right\}^{1/p}\leq CR^{\frac{d(d-1)(p-2)}{d(p-2)+p}}(\log R)^{\frac{(\gamma-2)(d-1)}{2(\gamma-1)}}.$$

At this point one can use interpolation to obtain a result for the remainder case. In particular, for $\gamma< d+1$, one has to interpolate between $p_{1}=2(\gamma-1)/(\gamma-2)$ and $\infty$.
We have 
\[
\|D_{R\Omega}(\sigma,t) \|_{\infty}\leq CR^{d(d-1)/(d+1)}
\]
 and 
\[
\|D_{R\Omega}(\sigma,t) \|_{p_{1}}\leq CR^{d(d-1)/(d+\gamma+1)}(\log R)^{\frac{(\gamma-2)(d-1)}{2(\gamma-1)}}.
\]
Hence 
\[
\|D_{R\Omega}(\sigma,t) \|_{\infty}^{(p-p_{1})/p} \|D_{R\Omega}(\sigma,t) \|_{p_{1}}^{p_{1}/p} \leq CR^{\frac{d(d-1)}{d+1}\left(1-\frac{2(\gamma-1)}{p(d+\gamma-1)}\right)}(\log R)^{\frac{d-1}{p}} .
\]

In the case $\gamma=d+1$, the critical indices for $p$ are all the same:
\[
\frac{2(\gamma-1)}{\gamma-2}=\frac{2d}{d-1}.
\]
The reader can easily resume the proof above to have the theorem.

Now we prove 2). First let $\gamma>d+1$ and $p<2(\gamma-1)/(\gamma-2)$. Because of (\ref{rot}), one has
\[
R^{d} \left\{ \sum_{n\neq 0} \left\| \widehat{\chi}_{\Omega}\right\|_{L^{p}(SO(d))}^{q}\right\}^{1/q}
 \leq R^{(d-1)/2}\left\{\sum_{n\neq 0} |n|^{-\frac{d+1}{2}q} \right\}^{1/q}\\
\leq R^{(d-1)/2}.
\]
Let $p=2(\gamma-1)/(\gamma-2)$. Because of (\ref{rot}), one has
\begin{align*}
& R^{d} \left\{ \sum_{n\neq 0} \left\| \widehat{\chi}_{\Omega}\right\|_{L^{p}(SO(d))}^{q}\right\}^{1/q}
\leq R^{(d-1)/2}\left\{\sum_{n\neq 0} |n|^{-\frac{d+1}{2}q}(\log(Rn))^{\frac{(\gamma-2)(d-1)}{2(\gamma-1)}q} \right\}^{1/q}\\
&\leq CR^{(d-1)/2}(\log R)^{\frac{(\gamma-2)(d-1)}{2(\gamma-1)}}\left\{\sum_{n\neq 0} |n|^{-\frac{d+1}{2}q} \right\}^{1/q} + 
CR^{(d-1)/2}\left\{\sum_{n\neq 0} |n|^{-\frac{d+1}{2}q}(\log n)^{\frac{(\gamma-2)(d-1)}{2(\gamma-1)}q} \right\}^{1/q}\\
&\leq CR^{(d-1)/2}(\log R)^{\frac{d-1}{p}}.
\end{align*}

Let $p>2(\gamma-1)/(\gamma-2)$. As before, using (\ref{rot}), if $p<\frac{(2d-1)\gamma-(d-1)}{(d-1)(\gamma-1)}$
\begin{align*}
R^{d} \left\{ \sum_{n\neq 0} \left\| \widehat{\chi}_{\Omega}\right\|_{L^{p}(SO(d))}^{q}\right\}^{1/q}&
 \leq R^{d-(d-1)\left(\frac{1}{p}+\frac{1}{\gamma}-\frac{1}{p\gamma} \right)-1}\left\{\sum_{n\neq 0} |n|^{-q(d-1)\left(\frac{1}{p}+\frac{1}{\gamma}-\frac{1}{p\gamma} \right)-q} \right\}^{1/q}\\
&\leq R^{(d-1)\left(1-\frac{1}{p}\right) \left(1-\frac{1}{\gamma} \right)}.
\end{align*}

Notice that $\frac{2(\gamma-1)}{\gamma-2}<\frac{(2d-1)\gamma-(d-1)}{(d-1)(\gamma-1)}$ if $\gamma>d+1$.

Let $p=\frac{(2d-1)\gamma-(d-1)}{(d-1)(\gamma-1)}$. One can use a cut-off function to give an estimate:
\begin{align*}
&R^{d}\left\{ \sum_{n\neq 0} (1+|\varepsilon n|)^{-K} \left\|\widehat{\chi}_{\Omega}(R\sigma(n)) \right\|^{q}_{L^{p}(SO(d))} \right\}^{1/q}\\
&\leq R^{d} \left\{\int_{|\xi|\geq 1} (1+\varepsilon |\xi|)^{-K}\left(R|\xi| \right)^{-(d-1)q\left(\frac{1}{p}+\frac{1}{\gamma}-\frac{1}{p\gamma} \right)-q}d\xi\right\}^{1/q}\\
&= R^{(d-1)\left(1-\frac{1}{p} \right)\left(1-\frac{1}{\gamma} \right)} \left\{\int_{1}^{+\infty} (1+\varepsilon \rho)^{-K}\rho^{-(d-1)q\left(\frac{1}{p}+\frac{1}{\gamma}-\frac{1}{p\gamma} \right)-q+d-1} d\rho\right\}^{1/q}\\
&= R^{(d-1)\left(1-\frac{1}{p} \right)\left(1-\frac{1}{\gamma} \right)}\varepsilon^{(d-1)\left(\frac{1}{p}+\frac{1}{\gamma}-\frac{1}{p\gamma} \right)+1-\frac{d}{q}} \left\{\int_{\varepsilon}^{+\infty} (1+y)^{-K}y^{-(d-1)q\left(\frac{1}{p}+\frac{1}{\gamma}-\frac{1}{p\gamma} \right)-q+d-1} dy\right\}^{1/q}\\
&\leq R^{(d-1)\left(1-\frac{1}{p} \right)\left(1-\frac{1}{\gamma} \right)}\varepsilon^{(d-1)\left(\frac{1}{p}+\frac{1}{\gamma}-\frac{1}{p\gamma} \right)+1-\frac{d}{q}}
\end{align*}
and with an appropriate $\varepsilon$ one obtains
\[
\left\{ \int_{SO(d)}\int_{T^{d}} |D_{R\Omega}(\sigma,t)|^{p}dt d\sigma\right\}^{1/p}\leq CR^{(d-1)\left(1-\frac{1}{p} \right)\left(1-\frac{1}{\gamma} \right)}.
\]

At this point one can use interpolation to obtain a result for the remainder case. In particular, for $\gamma > d+1$, we interpolate between $p_{2}=\frac{(2d-1)\gamma-(d-1)}{(d-1)(\gamma-1)}$ and $\infty$.

We have 
\[
\|D_{R\Omega}(\sigma,t) \|_{\infty}\leq CR^{(d-1)\left(1-\frac{1}{\gamma} \right)}
\]
 and 
\[
\|D_{R\Omega}(\sigma,t) \|_{p_{2}}\leq CR^{(d-1)\left(1-\frac{1}{\gamma}\right)\frac{d\gamma}{(2d-1)\gamma-(d-1)}}.
\]
Hence 
$$\|D_{R\Omega}(\sigma,t) \|_{\infty}^{(p-p_{2})/p} \|D_{R\Omega}(\sigma,t) \|_{p_{2}}^{p_{2}/p} \leq CR^{(d-1)\left(1-\frac{1}{p}\right) \left(1-\frac{1}{\gamma}\right)}.$$

\endproof

\proof[Proof of Corollary \ref{teo2}] For $1\leq p\leq 2$ it's true that $\|f\|_{p}\leq \|f\|_{2}=\|\widehat{f}\|_{2}$. Therefore 
\[
\left\{ \int_{SO(d)}  \int_{T^{d}} |D_{R\Omega}(\sigma,t)|^{p}dtd\sigma  \right\}^{1/p}\leq R^{d}\left\{\sum_{n\neq 0} \|\widehat{\chi}_{\Omega} \|^{2}_{L^{2}(SO(d))} \right\}^{1/2}.
\]
Because $2<2(\gamma-1)/(\gamma-2)$ for all $\gamma$, one has
\[
R^{d}\left\{\sum_{n\neq 0} \|\widehat{\chi}_{\Omega} \|^{2}_{L^{2}(SO(d))} \right\}^{1/2}\leq CR^{\frac{d-1}{2}}\left\{\sum_{n\neq 0} |n |^{-(d+1)} \right\}^{1/2}.
\]
Hence 
\[
\left\{ \int_{SO(d)}  \int_{T^{d}} |D_{R\Omega}(\sigma,t)|^{p}dtd\sigma  \right\}^{1/p}\leq CR^{\frac{d-1}{2}}.
\]
\endproof

\proof[Proof of Theorem \ref{teo3}] 
Using the mollified discrepancy with a cut-off function $\varphi$, one has
\[
|D_{R\Omega}(\sigma)|\leq CR^{d-1}\varepsilon+R^{d}\left\vert \sum_{n \neq 0} \widehat{\varphi}(\varepsilon \sigma(n))\widehat{\chi}_{\Omega}(R\sigma(n)) \right\vert.
\]
Therefore, by (\ref{rot})
\begin{align*}
\int_{SO(d)} R^{d}\left\vert \sum_{n \neq 0} \widehat{\varphi}(\varepsilon \sigma(n))\widehat{\chi}_{\Omega}(R\sigma(n)) \right\vert d\sigma &\leq
 R^{d}\sum_{n \neq 0} \widehat{\varphi}(\varepsilon \sigma(n)) \int_{SO(d)} |\widehat{\chi}_{\Omega}(R\sigma(n))|d\sigma\\
& \leq R^{\frac{d-1}{2}}\sum_{n \neq 0} (1+|\varepsilon n|)^{-K} |n|^{-\frac{d+1}{2}}.
\end{align*}
Using polar coordinates and then the change of variables $\varepsilon \rho=y$, one obtains
\begin{align*}
R^{\frac{d-1}{2}}\sum_{n \neq 0} (1+|\varepsilon n|)^{-K} |n|^{-\frac{d+1}{2}} & \leq R^{\frac{d-1}{2}} \int_{|\xi|\geq 1} (1+\varepsilon |\xi|)^{-K} |\xi|^{-\frac{d+1}{2}}d\xi\\
&=R^{\frac{d-1}{2}}\int_{1}^{+\infty} (1+\varepsilon \rho)^{-K} \rho^{-\frac{d+1}{2}+d-1}d\rho\\
&= R^{\frac{d-1}{2}}\varepsilon^{-\frac{d-1}{2}} \int_{\varepsilon}^{+\infty} (1+y)^{-K} y^{-\frac{d+1}{2}+d-1}dy
\end{align*}
Hence, because of the convergence of the last integral,
\[
\int_{SO(d)}|D_{R\Omega}(\sigma)|d\sigma \leq CR^{d-1}\varepsilon+R^{\frac{d-1}{2}}\varepsilon^{-\frac{d-1}{2}}
\] 
and with $\varepsilon=R^{-(d-1)/(d+1)}$ one has 
\[
\int_{SO(d)}|D_{R\Omega}(\sigma)| d\sigma \leq R^{\frac{d(d-1)}{d+1}}.
\]
\endproof

In order to estimate the discrepancy form below, we need first the following lemma.

\begin{lemma}\label{asin}
Let $\Omega \subset \mathbb{R}^{d}$ be a convex body with everywhere positive curvature except at one single point, which we may assume to be the origin, and with the boundary $\partial \Omega$ as the graph of the function $y=|x|^{\gamma}$ in a neighborhood of the origin, with $\gamma \geq 2$. Let $I$ be a small closed interval contained in $(0,\pi)$. For every direction $\theta \in I$, let $\sigma_{1}(\theta)$ and $\sigma_{2}(\theta)$ be the two points on $\partial \Omega$ where the tangents are perpendicular to $\Theta=(\cos\theta,\omega \sin\theta)$, with $\omega \in S^{d-2}$. We assume that the curvatures $K(\sigma_{1}(\theta))$ and $K(\sigma_{2}(\theta))$ are positive. Then
\begin{align*}
\widehat{\chi}_{\Omega}(\rho\Theta)&=\frac{1}{2\pi i}\rho^{-(d+1)/2}[K^{-1/2}(\sigma_{1}(\theta))\exp\left(-2\pi i\rho\Theta\cdot \sigma_{1}(\theta)-\pi i(d-1)/4\right) \\
& - K^{-1/2}(\sigma_{2}(\theta))\exp\left(-2\pi i\rho\Theta\cdot \sigma_{2}(\theta)+\pi i(d-1)/4\right) ] + O(\rho^{-(d+3)/2}).
\end{align*}
\end{lemma}

\proof See \cite{Herz1, Herz, Hlawka} or \cite[Corollary 7.7.15]{Hormander} or \cite[Chapter VIII]{Stein}. For a proof, see \cite[Lemma 2]{anelli}. \endproof

\proof[Proof of Theorem \ref{teo4}] 
Let $I\subset [\pi/2-\varepsilon, \pi/2+\varepsilon ]$ and let $\Theta=(\cos\theta,\omega \sin\theta)$, with $\omega \in S^{d-2}$ and $\theta \in I$. By our assumptions on $\Omega$, for all $\theta \in I$ there are two points $\sigma_{1}(\theta)$ and $\sigma_{2}(\theta)$ on $\partial \Omega$ where the tangents are perpendicular to $\Theta$. We can suppose that $\sigma_{1}(\theta)$ is the one close to the origin. In particular, we can say that there exists a positive constant $\kappa$ such that for all $\theta \in I$ is $K(\sigma_{2}(\theta))>\kappa$ and there exists an interval $J\subset I$ such that $K(\sigma_{1}(\theta))<\kappa/2$ for each $\theta \in J$, since $K(\sigma_{1}(\theta))\rightarrow 0$ as $\theta\rightarrow \pi/2$. Hence, by Lemma (\ref{asin})
\begin{align*}
&\int_{SO(d)} |\widehat{\chi}_{\Omega}(\rho \Theta)|d\Theta = \int_{SO(d-1)} \int_{0}^{\pi} |\widehat{\chi}_{\Omega}(\rho\cos\theta, \rho\omega\sin\theta)|\sin^{d-2}\theta d\theta d\omega\\
&> \int_{SO(d-1)} \int_{J \cup J+\pi}|\widehat{\chi}_{\Omega}(\rho\cos\theta, \rho\omega\sin\theta)|\sin^{d-2}\theta d\theta d\omega\\
&\geq c\rho^{(d+1)/2}\int_{SO(d-1)} \int_{J \cup J+\pi} \left(\left|K^{-1/2}(\sigma_{1}(\theta))-K^{-1/2}(\sigma_{2}(\theta)) \right|-c\rho^{-(d+3)/2}\right)\sin^{d-2}\theta d\theta d\omega\\
& \geq c\rho^{(d+1)/2}.
\end{align*}
Therefore, for all $k\neq 0$, one has 
\begin{align*}
\left\{\int_{SO(d)}\int_{T^{d}}|D_{R\Omega}(\sigma, t)|^{p}dt d\sigma\right\}^{1/p} &= R^{d}\left\{\int_{SO(d)}\int_{T^{d}} \left\vert \sum_{n\neq 0} \widehat{\chi}_{\Omega}(R\sigma(n))\exp(2\pi int)\right\vert^{p} dt d\sigma \right\}^{1/p}\\
&\geq R^{d} \left\{\int_{SO(d)} |\widehat{\chi}_{\Omega}(R\sigma(k))|^{p} d\sigma \right\}^{1/p} \geq R^{(d-1)/2}.
\end{align*}
\endproof

\end{document}